\documentclass[twoside,12pt]{amsart}
\usepackage{amsfonts}

\usepackage{graphicx}
\usepackage{amscd}
\usepackage{amsmath}
\newcommand{\chapter}{\section}

\begin{document}


\newtheorem{Thm}{Theorem}
\newtheorem{Ax}{Axiom}
\newtheorem{Prop}{Proposition}
\newtheorem{Cor}[Prop]{Corollary}
\newtheorem{Main}{}
\renewcommand{\theMain}{}
\newtheorem{Lem}[Prop]{Lemma}
\newtheorem{Fact}{Fact}
\renewcommand{\theFact}{}
\newtheorem{Claim}{Claim}
\renewcommand{\theClaim}{}

\newtheorem{Def}{Definition}
\newtheorem{rmk}{Remark}
\newenvironment{Rmk}{\begin{rmk}\em}{\end{rmk}}
\newtheorem{exm}{Example}
\newenvironment{Exm}{\begin{exm}\em}{\end{exm}}

\newtheorem{prf}{Proof}
\renewcommand{\theprf}{}
\newenvironment{Prf}{\begin{prf}\em}{\qed\end{prf}}
\newtheorem{prff}{}
\renewcommand{\theprff}{}
\newenvironment{Prff}{\begin{prff}\em}{\qed\end{prff}}



\newcommand{\YES}[1]{#1}
\newcommand{\NOT}[1]{}

\newcommand{\cA}{{\cal A}}
\newcommand{\cB}{{\cal B}}
\newcommand{\cC}{{\cal C}}
\newcommand{\cD}{{\cal D}}
\newcommand{\cE}{{\cal E}}
\newcommand{\cF}{{\cal F}}
\newcommand{\cG}{{\cal G}}
\newcommand{\cH}{{\cal H}}
\newcommand{\cI}{{\cal I}}
\newcommand{\cJ}{{\cal J}}
\newcommand{\cK}{{\cal K}}
\newcommand{\cL}{{\cal L}}
\newcommand{\cM}{{\cal M}}
\newcommand{\cN}{{\cal N}}
\newcommand{\cO}{{\cal O}}
\newcommand{\cP}{{\cal P}}
\newcommand{\cQ}{{\cal Q}}
\newcommand{\cR}{{\cal R}}
\newcommand{\cS}{{\cal S}}
\newcommand{\cT}{{\cal T}}
\newcommand{\cU}{{\cal U}}
\newcommand{\cV}{{\cal V}}
\newcommand{\cW}{{\cal W}}
\newcommand{\cX}{{\cal X}}
\newcommand{\cY}{{\cal Y}}
\newcommand{\cZ}{{\cal Z}}

\newcommand{\bbb}[1]{{\mathbb{#1}}}

\newcommand{\bN}{\bbb{N}}
\newcommand{\bZ}{\bbb{Z}}
\newcommand{\bR}{\bbb{R}}
\newcommand{\bC}{\bbb{C}}
\newcommand{\bQ}{\bbb{Q}}
\newcommand{\bT}{\bbb{T}}

\newcommand{\noind}[1]{{\setlength{\parindent}{0cm} #1}}
\newcommand{\parsk}{\par\medskip}
\newcommand{\pa}{\par\medskip}

\newcommand{\varend}{\end{document}}

\newcommand{\LP}{\left(}  \newcommand{\RP}{\right)}
\newcommand{\LQ}{\left[}  \newcommand{\RQ}{\right]}
\newcommand{\LB}{\left\{} \newcommand{\RB}{\right\}}
\newcommand{\LA}{\left\langle} \newcommand{\RA}{\right\rangle}
\newcommand{\LS}{\left|}  \newcommand{\RS}{\right|}
\newcommand{\LN}{\left\|} \newcommand{\RN}{\right\|}
\newcommand{\bLP}{\bigl(}  \newcommand{\bRP}{\bigr)}

\newcommand{\ts}{{\otimes}} \newcommand{\TS}{{\bigotimes}}

\newcommand{\beware}{~$\Psi\!\Psi\!\Psi$~}

\newcommand{\EM}{\em\bf}

\newcommand{\BE}{\begin{equation}}  \newcommand{\EE}{\end{equation}}
\newcommand{\BR}{\begin{eqnarray*}}  \newcommand{\ER}{\end{eqnarray*}}
\newcommand{\BER}{$$\begin{array}}  \newcommand{\EER}{\end{array}$$}
\newcommand{\head}[1]{
             \par\bigskip\noind{{\large\bf#1}\nopagebreak\par}\medskip}
\newcommand{\chead}[1]{\begin{center}
              \par\bigskip{\large\bf#1}\nopagebreak\par\medskip
              \end{center}}

\newcommand{\fillitem}{\L{\embox{}}\vspace{-.6cm}}

\newcommand{\toprove}[1]{{\buildrel\mbox{?}\over#1}}
\newcommand{\DEF}{{\buildrel\mbox{\scriptsize\,\,def\,\,}\over=}}
\newcommand{\RRe}{{\mbox{Re}\,}}
\newcommand{\IIm}{{\mbox{Im}\,}}
\newcommand{\BAR}[1]{\overline{#1}}
\newcommand{\HAT}[1]{{\buildrel{\,\,\wedge}\over{#1}}\!}
\newcommand{\TIL}[1]{{\buildrel{\,\sim\,}\over{#1}}}

\newcommand{\I}{\infty}
\newcommand{\es}{\emptyset}
\newcommand{\sm}{\setminus} \newcommand{\st}{\subset}
\newcommand{\eps}{\varepsilon}
\newcommand{\al}{\alpha} \newcommand{\be}{\beta}
\newcommand{\ga}{\gamma} \newcommand{\de}{\delta}
\newcommand{\DE}{\Delta}
\newcommand{\OM}{\Omega} \newcommand{\om}{\omega}
\newcommand{\la}{\lambda}
\newcommand{\half}{\frac{1}{2}}
\newcommand{\NI}{{n\to\infty}}
\newcommand{\ang}{{<\!\!\!)}}
\newcommand{\arc}[1]{\breve{#1}}

\newenvironment{txeqn}[1]{\begin{equation}\label{#1}\begin{array}{l}}{
                         \end{array}\end{equation}}

\setlength{\parindent}{0pt}






\title{An unusual proof that the reals are uncountable}
\author{Eliahu Levy}
\address{Department of Mathematics\\
Technion -- Israel Institute of Technology,
Haifa 32000, Israel}
\email{eliahu@techunix.technion.ac.il}


\date{}


\begin{abstract}
This somewhat unusual proof for the fact that the reals are uncountable,
which is adapted from one of Bourbaki's proofs in ``Fonctions d'une variable
r\'eelle'', may be of some interest.
\end{abstract}

\maketitle

\begin{Claim}
The set $\bR$ of the real numbers is uncountable.
\end{Claim}

\begin{Prf}
(adapted from a proof in \cite{Bourbaki}).
Suppose $\bR$ was countable.
Then there is a function $a(x):\bR\to\bR$ such that:

1. $a(x)>0$ for all $x$ (strictly positive),

2. the sum of the $a(x)$ on any finite set is $\le 1$.

(just take $a(x)=2^{-n}$ if $x$ is the $n$'th element in the counting).
\pa

Now, define for any subset $S$ of $\bR$,
\pa

$m(S)$ $:=$ the supremum of the sums of $a(x)$ on finite subsets of $S$.
\pa

Then surely $0\le m(S)\le 1$ for any $S$.
(In fact, $m(S)$ will be used only for $S$ an open ray.)
\pa

For $x\in\bR$, compare $m]-\infty,x[$ (open ray) with $x$. One may define:
$$c := \sup\Big\{x\,\Big|\, m]-\infty,x[\, > x\Big\},$$
and we shall reach a contradiction.
\pa

Indeed, since $a(c)>0$, there is a $y$ such that $y>c-a(c)$ and
$m]-\infty,y[ > y$, thus $y\le$ than the supremum $c$.
Now, since $y\le c$, $]-\infty,y[$ does not contain $c$.
But $]-\infty,y+a(c)[$ contains $c$ because $y>c-a(c)$.
So by the definition of $m(S)$,
$$m]-\infty,y+a(c)[\,\,\,\ge\,\,m]-\infty,y[\,+a(c) > y+a(c),$$
yet $y+a(c)$ exceeds the supremum $c$, a contradiction.
\end{Prf}

\end{document}